\documentclass[a4paper,CjK]{article}
\renewcommand{\normalsize}{\fontsize{10pt}{\baselineskip}\selectfont}
\usepackage[sf,sl,outermarks]{titlesec}
\titleformat{\section}{\large\bfseries}{\thesection}{1em}{}
\titleformat{\subsection}{\normalsize\bfseries}{\thesubsection}{1em}{}
\usepackage{ifpdf,appendix,caption}
\usepackage{amsmath,url}
\usepackage{amssymb}
\usepackage{enumerate}
\usepackage{amsthm,amsfonts}
\usepackage{graphicx}
\usepackage{color,verbatim}
\usepackage{mathrsfs,subeqnarray,subfigure}
\usepackage{listings}
\usepackage[ruled,vlined]{algorithm2e}
\newcommand{\nn}{\nonumber}
\newtheorem{theorem}{Theorem}[section]

\newtheorem{lemma}[theorem]{Lemma}

\newtheorem{proposition}[theorem]{Proposition}
\theoremstyle{definition}

\newtheorem{{defn}}[theorem]{Definition}

\graphicspath{{pic/}} 

\newcommand{\Rmn}[1]{\uppercase\expandafter{\romannumeral#1}}

\numberwithin{equation}{section}  
\numberwithin{table}{section}
\numberwithin{figure}{section}
\newcommand{\reff}[1]{(\ref{#1})}
\newcommand{\mA}{\mathcal{A}}
\newcommand{\mD}{\mathcal{D}}
\newcommand{\mL}{\mathcal{L}}
\newcommand{\mQ}{\mathcal{Q}}

\newcommand{\mbR}{\mathbb{R}}
\newcommand{\be}{\begin{equation}}
\newcommand{\ee}{\end{equation}}
\newcommand{\st}{\mathrm{s.t.}}

\newcommand{\eproof}{{\bf Q.E.D.}}
\def\lmd{{\lambda}}
\def\T{{\text{T}}}
\begin{document}

\title{Computing The Analytic Connectivity of A  Uniform Hypergraph}

\author{
  {Chunfeng Cui}\thanks{%
   Department of Electronic Engineering, College of Science and Engineering, City university of Hong Kong, Kowloon,
   Hong Kong ({\tt chunfcui@cityu.edu.hk}).
    } \and
  {Ziyan Luo}\thanks{%
  State Key Laboratory of Rail Traffic Control and Safety, Beijing
Jiaotong University, Beijing 100044, P.R. China; ({\tt
starkeynature@hotmail.com}).
   } \and
  {Liqun Qi}\thanks{%
    Department of Applied Mathematics, The Hong Kong Polytechnic University,
    Hung Hom, Kowloon, Hong Kong ({\tt maqilq@polyu.edu.hk}).
   }\and
      {Hong Yan}\thanks{%
  Department of Electronic Engineering, College of Science and Engineering, City University of Hong Kong, Kowloon, Hong Kong  ({\tt h.yan@cityu.edu.hk}).
    }
}

\maketitle

\begin{abstract} The analytic connectivity, proposed as a substitute of the algebraic connectivity in the setting of hypergraphs, 
is an important quantity in spectral hypergraph theory. The definition of the analytic connectivity for a uniform hypergraph 
involves a series of polynomial
 optimization problems (POPs) associated with the Laplacian tensor of the hypergraph   with  nonnegativity 
 constraints and a sphere constraint, which poses  difficulties in computation. To reduce the involved computation, 
 properties on the algebraic connectivity are further exploited, and several important structured uniform hypergraphs 
 are shown to attain their analytic connectivities at vertices of the minimum degrees,  hence admit a relatively 
 less computation by solving a small number of POPs. To efficiently solve each involved POP, we propose a feasible 
 trust region algorithm ({\tt FTR}) by exploiting their special structures. The global convergence of {\tt FTR} 
 to the  second-order necessary conditions points is established,  and numerical results for both small and large 
 size examples with comparison to other existing algorithms for POPs are reported to demonstrate the efficiency 
 of our proposed algorithm. 

\vskip 12pt \noindent {\bf Key words.} {Uniform hypergraph; Laplacian tensor; Analytic connectivity; Feasible trust region algorithm}

\vskip 12pt\noindent {\bf AMS subject classifications. }{05C65, 15A18, 90C55}
\end{abstract}
\section{Introductions}

Spectral graph theory is a well-studied and highly applicable
subject, which focuses on the connection between properties of a graph and the eigenvalues
of matrices associated with the graph.   Such matrices include the adjacency matrix, the Laplacian matrix and the signless Lpalacian matrix of the graph \cite{BM, Ch, CDGT, CDS, LSG}. However, the study of graphs cannot fully meet the developments of modern science and technology, especially in big data analysis and complex networks.   This motivates the study of hypergraphs, where  an edge may connect more than two vertices \cite{Be, Br}, comparing to two-vertices edges in ordinary graphs.  Spectral hypergraph theory correspondingly emerged which was based upon matrix spectral analysis in its early stage.
\vskip 2mm

In 2005, Lim \cite{Li} and Qi \cite{Qi0} independently introduced the concept of eigenvalues for tensors, which initiated the study of tensor spectral theory and   paved a way for the development of spectral hypergraph theory via tensors.
The related research include
spectral hypergraph theory \cite{CD,CD1,GD1,KF,  KFT, LSQ,Pe,PZ1,SQH,SSW,XQ1,ZSWB},
eigenvalues \cite{HQ1,LQY,PZ,XC,XC1,XC2,XQ,YQS,YZL},
connectivity \cite{HQ, CCL},
Laplacian tensor \cite{BFZ,HQ2,HQX,PZ2,Qi, YZL1},
 structured tensors related \cite{CQ,CLN},
 special hypergraphs \cite{BZW,HQS,KNY,QSW,YSS},
 hypergraph properties \cite{BP,FTPL,GD,Ni,Ni1}.
The tensors studied in these papers include adjacency tensors, Laplacian tensors and signless Laplacian tensors of hypergraphs. Benefitting from  the high sparsity of these tensors,  Chang, Chen and Qi \cite{CCQ}   recently proposed a CEST algorithms for computing extremal eigenvalues of large-scale adjacency tensors, Laplacian tensors and signless Laplacian tensors of uniform hypergraphs, which provides a  useful computational tool for spectral hypergraph theory via tensors.

\vskip 2mm

It is well-known that in spectral graph theory, the algebraic connectivity \cite{Fi}, defined as the second smallest eigenvalue of the Laplacian matrix of a graph, is an important quantity. However, as Laplacian tensors of uniform hypergraphs may have complex eigenvalues, a different approach for generalizing this concept to hypergraphs was introduced by Qi in \cite{Qi}, where the analytic connectivity for a uniform hypergraph was defined via an optimization formulation
\begin{align*}
  \alpha(G)=\min_{j=1,\cdots,n} \min_{x\in\mbR^n}\{\mathcal{L}x^k\,:\,x\in \mathbb{R}^n_+,\ \sum_{i=1}^nx[i]^k=1,\ x[j]=0\}.
\end{align*}
This is  shown to be linked with the edge connectivity of the hypergraph. It was further studied by Li, Cooper and  Chang in \cite{CCL} where the analytic connectivity was shown to be connected with   other important invariants of hypergraphs, such as the degree, the vertex connectivity, the diameter and the isoperimetric number.

\vskip 2mm

To our best knowledge, no efficient algorithm has been proposed for computing
the analytic connectivity of a uniform hypergraph in the literature.
The definition of analytic connectivity involves a series of polynomial optimization problems.
For   dimension $n$ large enough, this is very costly.
To fix this issue, we firstly explore  some specific hypergraphs,
and shown some properties on the vertices the analytic connectivity will possibly attained.

\vskip 2mm

For each specific $j$, we propose a feasible trust region algorithm  for the computation of this quantity.
Note that the analytic connectivity involves a series of optimization problems, each of which possesses nonnegativity constraints and a sphere  constraint. Thus, they are special cases of the following general constrained optimization problem
\begin{equation}\label{equ:gnp}
  \min_{x\in \mathbb{R}^n} \quad f(x)\quad \st\quad c(x)=0,\quad x\ge0,
\end{equation}
where $f:\mathbb{R}^n\rightarrow \mathbb{R}$ is a nonconvex polynomial function and  $c:\mathbb{R}^n\rightarrow \mathbb{R}$ is a nonlinear smooth function.
Existing optimization algorithms for \reff{equ:gnp} can be roughly classified into three types.
The first type includes the penalty methods, which incorporates the equality constrains $c(x)=0$ into the objective function as a penalty term,
 and attempt  to solve \reff{equ:gnp} by a sequential minimization problems of the form
\begin{equation*}
  \min_{x\in \mathbb{R}^n}\quad \mathcal{F}^{(1)}_k(x,\lambda_k) \quad \st\quad x\ge0,
\end{equation*}
where the objective function $\mathcal{F}^{(1)}_k$  could be any penalty function such that the subproblem can be easily solved. For instance, the L-BFGS method \cite{LBFGS-B}, the gradient projection method \cite{CM1987, DF2006}, and the active set method \cite{HZ2006}, etc. The solver {\tt MINOS} belongs to this type. However, as the hard constraint $c(x)=0$ has been relaxed as a penalty term in the objective, this type of methods usually result in an infeasible point for our sphere constraint. The second type of methods involves solving
\begin{equation*}
  \min_{x\in \mathbb{R}^n}\quad \mathcal{F}^{(2)}_k(x,\mu_k) \quad \st\quad c(x)=0,
\end{equation*}
where the objective function $\mathcal{F}^{(2)}_k$ is always with some interior-point penalty of the nonnegativity constraint and the solver {\tt IPOPT} belongs to this type. With the equality constraint in the above subproblem, this type of methods is always time consuming.
The third type includes the sequential quadratic programming methods, which solves the subproblem
\begin{equation*}
  \min_{x\in \mathbb{R}^n} \quad  \mathcal{F}^{(3)}_k(x,\lambda_k) \quad \st \quad  c_k+A_k(x-x_k)=0,\quad x\ge0,
\end{equation*}
where the objective function $\mathcal{F}^{(3)}_k$ is a quadratic function using the information of the Lagrangian function or its variants \cite{FS2005}. This type of methods show their strength when the constraints have significant nonlinearity, and the solver {\tt SNOPT} belongs to this type.

\vskip 2mm

Note that the  constraint $c(x)=0$ in this paper is actually the $k$-norm sphere constraint.
By exploring this special structure, we propose a feasible trust region method ({\tt FTR}), the mixture of trust region method and the projection method, in which the projection step ensures the feasibility of each iteration and the trust region technique enhances the convergence.
{\tt FTR} was  also used in \cite{hao2015feasible} for computing  Z-eigenvalues of symmetric tensors. While the main difference is that here we adopted the $\infty$-norm trust region instead of the  Euclidean norm, 
which remarkably facilitates the computation as at each iteration only a linear constrained quadratic subproblem needs to be handled.
Infinity norm was also used in \cite{GMT} for bound constrained problems, where advantages in terms of computational costs were demonstrated.


\vskip 2mm

This paper is organized as follows. In Section \ref{section:pre}, several related basic concepts and properties on hypergraphs and the analytic connectivity are reviewed. Further properties on the vertices attainable for the analytic connectivity is discussed in Section 3 to reduce the computation by cutting down the number of the involved POPs. For each POP, an {\tt FTR}  algorithm for computing the analytic connectivity of a uniform hypergraph is proposed in Section \ref{section:FTRalg}. The global convergence to the second order stationary points is established in Section \ref{section:convergence}. Numerical results are reported in Section \ref{section:numericalexperiments}, which demonstrates the efficiency of our proposed algorithm, and indicates that the analytic connectivity is a good choice to characterize the connectivity of the involved hypergraph as well. Conclusions are drawn in Section 7.

\vskip 2mm

Notations throughout the paper are listed here. Let $k$ and $n$ be any two positive integers. We use $T_{k,n}$ to denote the space of all $k$-th order $n$-dimensional tensors. $\mathbb{R}^n_+$ is used to stand for the set of all nonnegative vectors in $\mathbb{R}^n$. For any $x\in\mathbb{R}^n$ and any integer $i\in [n]$, $x[i]$ denotes the $i$-th  component of $x$, $x^{[k]}:=(x[i]^k)\in\mathbb{R}^n$ with any given positive integer $k$, and $\text{diag}(x)\in \mathbb{R}^{n\times n}$ is the diagonal matrix generated by $x$. For any set $C$, $|C|$ denotes the cardinality of $C$. The index set $\{1,2,\ldots,n\}$ is simply denoted as $[n]$. The notation $\binom{n}{m}$ denotes the combinatorial number of choosing $m$ from $n$.


\section{Preliminaries}
\label{section:pre}
As a natural extension of a graph, a {\bf hypergraph} $G = (V,E)$ with the vertex set $V=[n]$ and the edge set $E=\{e_1,\cdots,e_m\}$ allows each of its edge $e_j$ joins any number of vertices.  If  each edge $e_j$ connects  exactly $k$ vertices, this hypergraph is called a {\bf $k$-uniform hypergraph}, or simply called as a {\bf $k$-graph}. For more details on hypergraphs,   refer to \cite{Be, Br, CD}. Obviously, $G$ is reduced to an ordinary graph when $k=2$. Thus, we assume $k\ge3$ throughout the paper.
\vskip 2mm

Many important structured hypergraphs have been introduced in the literature. Let $G=(V,E)$ be a uniform hypergraph. $G$ is called a {\bf sunflower} if there is a disjoint partition of the vertex set $V$ as $V=V_0\cup V_1\cup \cdots \cup V_d$ such that $|V_0|=1$ and $|V_1|=\cdots=|V_d|=k-1$, and $E=\{V_0\cup V_i\,|\, i\in[d]\}$ (\cite{HQS}); $G$ is called a {\bf hypercycle} if there are $s$ subsets $V_1$, $\ldots$, $V_s$ of the vertex set $V$ such that $|V_1| = \cdots = |V_s| = k$, $|V_1 \cap V_2|=$ $\cdots$ $=|V_{s-1} \cap V_s|=$ $|V_s \cap V_1|=1$ and $V_i\cap V_j = \emptyset$ for the other cases, the intersections $V_1 \cap V_2$, $\ldots$, $V_s \cap V_1$ are mutually different, and $E = \{V_i \mid i \in [s]\}$ (\cite{HQX}); $G$ is called a {\bf squid} if we can number the vertex set $V$ as $V=\{i_{1,1},\cdots, i_{1,k}, \cdots, i_{k-1,1},\cdots,i_{k-1,k},i_{k,1}\}$ such that the edge set $E=\{\{i_{1,1},\cdots, i_{1,k}\}$, $\cdots$,$\{ i_{k-1,1},\cdots,$ $i_{k-1,k}\}$, $\{i_{1,1}, \cdots, i_{k-1,1}, i_{k,1}\}\}$ (\cite{HQS}); More generally, $G$ is called a {\bf $s$-path} of length $l$ if $V=\{v_1,v_2,\ldots, v_{s+l(k-s)}\}$ and $E=\{\{v_{1+i(k-s)},v_{1+i(k-s)+1},\ldots,v_{s+(i+1)(k-s)}\}\mid 0\leq i\leq l-1\}$. Particularly, we call a $1$-path hypergraph $G$ as a {\bf loose path}; $G$ is called a {\bf complete $k$-graph} if $E=\{e\mid e\subset V, |e|=k\}$.

\vskip 2mm

Some related fundamental concepts of uniform hypergraphs are reviewed as follows.


   \begin{{defn}}[\cite{CD, Qi}] Let $G=(V, E)$ be a $k$-graph.
   The \textbf{adjacency tensor} of $G$ is
   defined as the $k$-th order $n$-dimensional tensor $\mathcal{A}$ whose $(i_1,\cdots,i_k)$-entry is:
   \begin{equation*}
        a_{i_1,\cdots,i_k}:=\left\{
                              \begin{array}{cl}
                                \frac{1}{(k-1)!}, & \hbox{\text{if}\ $\{i_1\cdots,i_k\}\in E;$} \\
                                0, & \hbox{\text{otherwise}.}
                              \end{array}
                            \right.
   \end{equation*}
   Let $\mathcal{D}$ be a $k$-th order $n$-dimensional diagonal tensor with its diagonal element $d_{i,\cdots,i}$
   being $d_i$, the degree of vertex $i$, for all $i\in[n]$. $\mD$ is called the \textbf{degree tensor}
   corresponding to $G$.
   Then \textbf{Laplacian tensor}  of $G$ is defined as $\mL: = \mD - \mA$, and the \textbf{signless Laplacian tensor} of $G$ as $\mQ:=\mD+\mA$.
   \end{{defn}}


\begin{{defn}}[\cite{Qi}]\label{Def:analytic-connectivity}
  Let $G$ be a $k$-graph with $n$ vertices. The {\bf analytic connectivity} of $G$  is defined as
\begin{align}\label{Def:anacon}
  \alpha(G)=\min_{j=1,\cdots,n}\alpha_j(G),
\end{align} where
 \begin{align}\label{Def:anaconj}
  \alpha_j(G)=\min_{x\in\mbR^n}\{\mathcal{L}x^k\,:\,x\in \mathbb{R}^n_+,\ \sum_{i=1}^nx[i]^k=1,\ x[j]=0\},
\end{align}
with $\mathcal{L}$ the Laplacian tensor of $G$.
\end{{defn}}


Let $G = (V, E)$ be a $k$-graph with $n$ vertices. For each vertex $i\in V$, denote by $E(i)$ the set of edges containing
the vertex $i$, i.e., $E(i):=\{e\in E\, | \, i\in e\}$. The degree $d_i$ of the vertex $i$ is the cardinality $|E(i)|$
of the set $E(i)$. Denote by $\Delta$, $\delta$ and $\bar{d}$ the maximum, minimum and average degree of $G$, respectively.
Existing results on analytic connectivity of a uniform hypergraph include the following:
 \begin{itemize}
   \item \cite{Qi} $\alpha(G)\ge0$; $\alpha(G)>0$ if and only if $G$ is connected;
   \item \cite{Qi} $e(G)\ge\frac{n}{k}\alpha(G)$, where $e(G)$ is the {\bf edge connectivity} of $G$, defined as the minimum cardinality of an edge cut of $G$;
   \item \cite{Qi} $\alpha(G)\le \delta$;
   \item \cite{CCL} $\alpha(K_n^{(k)}) = \left({n-2 \atop k-2}\right)$, where $K_n^{(k)}$ is the complete $k$-graph;
   \item \cite{CCL} denote $v(G)$ as the {\bf vertex connectivity} of $G$, defined as the minimum cardinality of a vertex cut of $G$,
   \begin{equation}\alpha(G) \le \bar{\alpha}:=\left({n-2 \atop k-2}\right)- \left[\left({n-v(G)-1 \atop k-1}\right) - \left({\lfloor{(n-v(G) \over 2}}\rfloor-1 \atop k-1\right)\right]{k-1 \over n-1};\label{upperbound}\end{equation}
     \item \cite{CCL} ${k \over 2}i(G) \ge \alpha(G) \ge \Delta - \sqrt{\Delta^2 - i^2(G)},$ where $i(G)$ is the {\bf isoperimetric number}, or the {\bf Cheeger constant} of $G$, defined by $i(G) = \min \left\{ {|E(S, \bar S)| \over |S|} : S \subset V, 0 < |S| \le {n \over 2} \right\},$ $\bar S = V \setminus S$, and $E(S, \bar S)$ is an edge cut of $G$;
   \item \cite{CCL} $\alpha(G)\ge \frac{4}{n^2(k-1){\rm diam}(G)}$, where diam$(G)$ is the {\bf diameter} of $G$, defined as the maximum distance between any pair of vertices of $G$;
   \item \cite{CCL} $\alpha(G)\le \min\left\{ \frac1k (d(v_{i_1}) + d(v_{i_2})+\cdots +d(v_{i_k})-k):
   {v_{i_1},\cdots,v_{i_k}\in E(G)}\right\}.$
 \end{itemize}

It is worth pointing out that the isoperimetric number or the Cheeger constant of an ordinary graph provides a numerical measure of whether or not a graph has a ``bottleneck", which has wide applications such as in constructing well-connected networks of computers and card shuffling. However, the computation of such an invariant is very difficult and the algebraic connectivity provides a reasonable good bound in terms of the well-known ``{\bf Cheeger inequality}" in the ordinary graph case. This result is in a certain sense theoretically extended to the uniform hypergraphs as stated above by  Li, Cooper and  Chang \cite{CCL} where the analytical connectivity was adopted instead of the algebraic connectivity. In this regard, the computational algorithm presented in this paper makes the theoretical result of \cite{CCL} practically feasible to efficiently bound the isoperimetric number of a $k$-graph.

\section{Properties on the analytic connectivity}
\label{section:Property}
In this section, we will discuss the properties on finding which vertices of a uniform hypergraph the analytic connectivity will possibly be attained at. This will henceforth play an essential role in reducing the required computation for the analytic connectivity by cutting down the number of POPs involved in Definition \ref{Def:analytic-connectivity}. We begin with the following important lemma.

 \begin{lemma}\label{lem:monotone}
  Let $G=(V,E)$ be a $k$-graph with $V=[n]$, and $i, j\in [n]$ be any two vertices with edge sets $E(i)$ and $E(j)$. If $E(i)\subset E(j)$, then $
  \alpha_i(G) \leq \alpha_j(G)$, where $\alpha_i(G)$ and $\alpha_j(G)$ are defined as in (\ref{Def:anaconj}).
  \end{lemma}
  \textbf{Proof:}\quad
Let $E(i)=\{e_1(i),\ldots,e_{d_i}(i)\}$ and $E(j)=E(i)\cup \{e_{d_i+1}(j),\ldots,e_{d_j}(j)\}$, where $d_i$ and $d_j$ are the degrees of vertices $i$ and $j$, respectively. For any $x\in\mbR^n$, denote $\mathcal{L}(e)x^k = \sum_{i\in e}x[i]^k-k\Pi_{i\in e}x[i]$ as the Laplacian function corresponding to any given edge $e\in E$. For any $x_1\in\mbR^n$ satisfying $x_1[i]=0$, we have
 $$\mathcal{L}x_1^k=\sum_{e_l\in E(i)}\sum_{l_t\in e_l, l_t\neq i} x_1[l_t]^k +
 \sum_{e\in E \setminus E(i)}\mathcal{L}(e)x_1^k. $$
For any $x_2\in\mbR^n$ satisfying $x_2[j]=0$, we have
 $$\mathcal{L}x_2^k=\sum_{e_l\in E(i)}\sum_{l_t\in e_l, l_t\neq j} x_2[l_t]^k +
 \sum_{e_l\in E(j)\setminus E(i)}\sum_{l_t\in e_l, l_t\neq j} x_2[l_t]^k +
 \sum_{e\in E \setminus E(j)}\mathcal{L}(e)x_2^k. $$
Note that the vertex $i$ is only contained in the edges of $E(i)$ and hence $x_2[i]$  only exists in the first term of the right hand side of the above expression. To achieve the minimum value $\alpha_j(G)$ in \reff{Def:anaconj}, it is evident from the nonnegativity constraint that for any optimal solution $\bar{x}$ of the problem
\reff{Def:anaconj} with $x[j]=0$, it holds that $\bar{x}[i]=0$. Therefore, $\bar{x}$ is also a feasible solution of the problem
\reff{Def:anaconj} with $x[i]=0$. This immediately shows the desired inequality.
  \quad\eproof


\vskip 2mm
With the help of Lemma \ref{lem:monotone}, we can show that for several important uniform hypergraphs,
such as sunflowers, hypercycles,  squids and loose path,
the computation of their analytic connectivities can be significantly reduced by the following theorem.

\begin{theorem}\label{anacon:hyperstarcirclesun}
Let $G$ be a $k$-graph with the vertex set $[n]$. If $G$ is a sunflower, or a hypercycle, or a squid, or a loose path, then $\alpha(G)=\alpha_j(G)$, where $j\in [n]$ is a vertex with the minimum degree.
\end{theorem}
    \textbf{Proof:} \quad
   Let $G=(V,E)$ be a $k$-graph with $V=[n]$. (i) If $G$ is a sunflower, then we can find a disjoint partition of the vertex set $V$, says $V=V_0\cup V_1\cup \cdots \cup V_d$, such that $|V_0|=1$ and $|V_1|=\cdots=|V_d|=k-1$, and $E=\{V_0\cup V_i\,|\, i\in[d]\}$, where $1+d(k-1)=n$. Let $V_0=\{v_0\}$. Obviously,  $v_0$ has degree $d$ and other vertices all have degree $1$. Moreover, for any $v\in V\setminus V_0$,
   $E(v)\subset E(v_0)$. Invoking of Lemma \ref{lem:monotone}, the desired result follows readily in this case. (ii) If $G$ is a hypercycle, then
   there exist $s$ subsets $V_1$, $\ldots$, $V_s$ of the vertex set $V$ such that $|V_1| = \cdots = |V_s| = k$, $|V_1 \cap V_2|=$ $\cdots$ $=|V_{s-1} \cap V_s|=$ $|V_s \cap V_1|=1$ and $V_i\cap V_j = \emptyset$ for the other cases. From the definition of hypercycles, we know that each intersected vertex has degree two and others has degree one. And for any $v\in V$ of degree two, there exists a vertex $v'\in V$ such that $E(v')\subset E(v)$. Thus, by applying Lemma \ref{lem:monotone}, the desired result is obtained in this case. (iii) If $G$ is a squid, then we can number the vertex set $V$ as $V=\{i_{1,1},\cdots, i_{1,k}, \cdots, i_{k-1,1},\cdots,i_{k-1,k},i_{k,1}\}$ such that the edge set $E=\{\{i_{1,1},\cdots, i_{1,k}\}$, $\cdots$,$\{ i_{k-1,1},\cdots,$ $i_{k-1,k}\}$, $\{i_{1,1}, \cdots, i_{k-1,1}, i_{k,1}\}\}$.  Note that the vertices $i_{1,1}, \cdots, i_{k-1,1}$ all have degree two and others all have degree one, and for every vertex $i_{j,1}$ with degree two,
    there exist vertex $i_{j,2}$ such that $E(i_{j,2}) \subset E(i_{j,1})$. Thus, from Lemma \ref{lem:monotone}, we have $\alpha_{j,2}(G) \le \alpha_{j,1}(G)$. (iv) Similar to case (ii), we can prove the case when $G$ is a loose path by definition and Lemma \ref{lem:monotone}. This completes the proof.
    \quad\eproof

\vskip 2mm
Two more specific uniform hypergraphs are discussed whose analytic connectivities can be computed via solving \reff{Def:anaconj} with special choices of $j$. The first one is the $2$-path with $n$ vertices which is plotted as follows.

 \begin{figure}[h]
  \centering
  \includegraphics[width=2.5in]{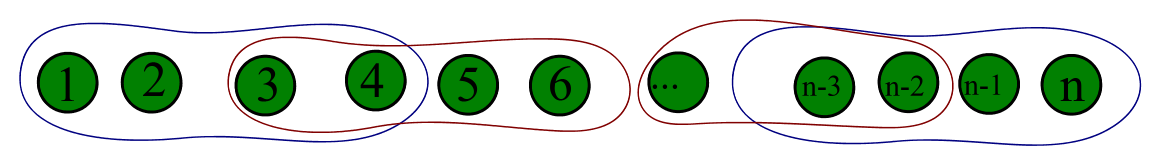}
  \caption{A $2$-path $4$-graph with length $\frac{n-2}{2}$}\label{fig:example10}
\end{figure}

\begin{proposition}\label{lemma:chaindecrease} Let $G$ be a $2$-path $4$-graph with $n\geq 4$ vertices, defined as in Figure \ref{fig:example10}. Then $\alpha(G)=\alpha_j(G)$, where $j$ could be any element in $\{1, 2, n-1, n\}$. Moreover, $\alpha(G)$ is monotonically decreasing with $n$. 
\end{proposition}
\textbf{Proof:\quad}
First we consider the first part of the proposition. It is trivial when $n=4$. For $n=6$, the desired result can be obtained immediately from the symmetric structure of $G$ and Lemma \ref{lem:monotone}.
Before proceeding for general cases of $n>6$, we will introduce the following useful function for any given even integer $l\ge 4$,
 \begin{align*}
 \beta_{\gamma}^{l} = \min_y\ \{g_{\gamma}^{l}(y):=
  y_1^4+\cdots +y_{l-2}^4-4y_1y_2y_3y_4-\cdots-4y_{l-3}y_{l-2}y_{l-1}y_{l}\quad
  \st\quad\sum_{i=1}^ly_i^4=\gamma\}.
 \end{align*}
It is easy to see that $\beta_{\gamma}^{l} = \gamma\beta_{1}^{l}$ from the homogeneous structure of the above minimization problem. Moreover, we claim that $\beta_{\gamma}^{l}$ is decreasing with $l$. Let $l_1$, $l_2$ be any two even integers satisfying $l_1>l_2\geq 4$. For any optimal solution $\hat{y}$ of problem with dimension  $l_2$, $\bar{y}=[zeros(l_1-l_2,1),\hat{y}]$ is a feasible solution of dimension $l_1$. Hence
 $$\beta_{\gamma}^{l_1}\le g_{\gamma}^{l_1}(\bar{y})=g_{\gamma}^{l_2}(\hat{y})=\beta_{\gamma}^{l_2},$$
 where the first equality comes from the fact that the formulation of $g_{\gamma}^{l_2}(\hat{y})$ is the
 same with $g_{\gamma}^{l_1}(\bar{y})$. Furthermore, for any even integer $l\geq 4$, $\beta_{\gamma}^{l}$ is negative. This comes from the claim above and the observation that
 given $\bar{y}=(\frac{\sqrt{2}}{2},\frac{\sqrt{2}}{2},\frac{\sqrt{2}}{2},\frac{\sqrt{2}}{2})$,
 $\beta_{\gamma}^{4}\le g_{\gamma}^{4}(\bar{y})=\bar{y}_1^4+\bar{y}_2^4-4\bar{y}_1\bar{y}_2\bar{y}_3\bar{y}_4=-\frac{1}{2}$.

 For any  even integer $n>6$, it holds that $\mathcal{L} x^4  = 1 + \sum_{i=3}^{n-2} x_i^4  - \mathcal{A} x^4.$
 Suppose that for some $j$, $x_j=0$, then the index set $[n]\setminus \{j\}$ can be partitioned into $\{1,\cdots,L_1\}$, $\{j-(-1)^j\}$, and $\{L_1+3,\cdots,n\}$. Hence $\mathcal{L} x^4$ can be rewritten as
 \begin{align}\label{equ:departchain}
 \mathcal{L} x^4 = 1+g_{\gamma}^{L_1}(x_{[1:L_1]})+\delta+g_{1-\gamma-\delta}^{n-L_1-2}(x_{[L_1+3:n]}),
 \end{align}
 where $L_1=j-\frac32-\frac12(-1)^j$,   $\delta=x_{j-(-1)^j}^4$, and $g_{\gamma}^2=0$.
 Note that the variable $x$ in \reff{equ:departchain} are partitioned into three subvectors, thus
 \begin{align*}
   \min\  \mathcal{L} x^4 \ \st \sum_{i=1, i\neq j}^n x_i^4=1 \Longleftrightarrow \min_{\gamma,\delta\ge0}\  1+\gamma \beta_{1}^{L_1}+\delta+(1-\delta-\gamma)\beta_1^{n-L_1-2}\ \st\  \gamma+\delta\le 1.
 \end{align*}
It follows from  t $\beta_{1}^{l}$ is negative that $\delta=0$,
and the objective function is reduced to $1+\gamma\beta_1^{L_1}+(1-\gamma)\beta_1^{n-L_1-2}$,
 as $\beta_1^l$ decreasing with $l$, hence
\begin{align*}
   \alpha_j(G)= 1+\beta_1^{l_j},
\end{align*}
where $l_j=\max(L_1, n-L_1-2)$.
Hence $ j^*= \arg\min_j \ \alpha_j = \arg\max_j\ l_j. $
By direct computation we have $j^*\in\{1,2,n-1,n\}$ and  $l_j =n-2$.
When $l_j=L_1$, it holds that $\gamma =1$;
otherwise, $\gamma =0$. Thus,
\begin{align}\label{equ:alhavalue}
  \alpha(G) = 1+\beta_1^{n-2}.
\end{align} As $\beta_1^{n-2}$ is monotonically decreasing with $n$,
so is the analytic connectivity $\alpha$ from (\ref{equ:alhavalue}). This completes the proof.
\eproof

%


\vskip 2mm

\begin{figure}[h]
  \centering
  \includegraphics[width=1.5in]{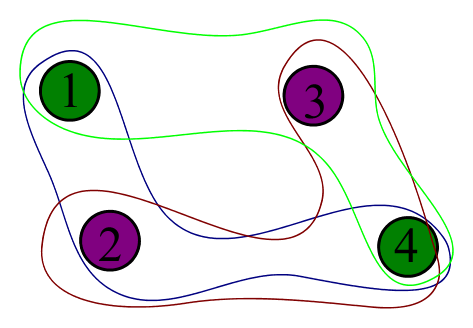}
  \caption{$K_4^{-}$ generated by deleting the edge $\{1, 2, 3\}$ from $K_4^{(3)}$} 
  \label{fig:examplecomplete}
\end{figure}

The second specific one, termed as $K_n^{-}$, is the $k$-graph obtained by deleting an arbitrary edge from a complete $k$-graph $K_n^{(k)}$.
For example, when $k=3$, $n=4$, the edge set of $K_4^{-}$ are $\{\{1, 2, 4\}$, $\{1, 3, 4\}$, $\{2, 3, 4\}\}$,
as shown in Figure \ref{fig:examplecomplete}.

\begin{proposition}\label{lem:kn-}
  Suppose $K_n^-$ is the hypergraph generated by deleting an edge $\hat{e}$ from $K_n^{(k)}$.
  Then $\alpha(K_n^-)=\alpha_j(K_n^-)$, where $j$ is some vertex in $\hat{e}$, i.e.,
  \begin{align}\label{equ:j1lj2}
    \alpha_{j_1}(K_n^-) < \alpha_{j_2}(K_n^-), ~~\forall j_1\in e,~ \forall j_2\in V\setminus \hat{e}.
  \end{align}
\end{proposition}
  {\bf Proof: \quad} Without loss of generality, suppose that the edge $\hat{e}=\{1,\cdots,k\}$ is deleted.
   By the symmetric property of this hypergraph, to show \reff{equ:j1lj2}, we only need to prove $\alpha_{1}(K_n^-) < \alpha_{n}(K_n^-)$.
  For $j\in[n]$ satisfying $x[j]=0$ and $\|x\|_k=1$, we have
    \begin{align*}
     \mathcal{L}(K_n^-) x^k&= \sum_{e\notin E(j)}\mathcal{L}(e) x^k+ \sum_{e\in E(j)}\mathcal{L}(e) x^k.
     \end{align*}
    If $j\notin\hat{e}$,
         \begin{align}\label{equ:jnotinE}
   \sum_{e\in E(j)}\mathcal{L}(e) x^k= \sum_{e\in E(j)}\sum_{l\in e, l\neq j}x[l]^k= \binom{n-2}{k-2}\sum_{l\in [n]\setminus j} x[l]^k = \binom{n-2}{k-2},
   \end{align}
   otherwise,
   \begin{align}\label{equ:jinE}
    \sum_{e\in E(j)}\mathcal{L}(e) x^k= \sum_{e\in E(j)}\sum_{l\in e, l\neq j}x[l]^k- \mathcal{L}(\hat{e})x^k= \binom{n-2}{k-2}-\mathcal{L}(\hat{e})x^k. 
   \end{align}
    
    For the case   $x[1]=0$, set $\bar{x}$ as
    \begin{align*}
        \bar{x}[i]=\left\{ \begin{array}{cl}
          (\frac{1}{n-1})^{\frac1k},&\quad \text{if} \quad i=2,\cdots,n;\\
          0,                         &\quad \text{if} \quad i=1.
        \end{array}
        \right.
    \end{align*}
Then $j\in\hat{e}$. For all $e\notin E(j)$ it holds that $\mathcal{L}(e) \bar{x}^k=0$.  It follows from \reff{equ:jinE} that
    \begin{align*}
      \alpha_1(K_n^-) \le \mathcal{L}(K_n^-) \bar{x}^k=
     \sum_{e\notin E(j)}\mathcal{L}(e) \bar{x}^k+ \sum_{e\in E(j)}\mathcal{L}(e) \bar{x}^k - (\bar{x}[2]^k+\cdots+\bar{x}[k]^k)
      = \binom{n-2}{k-2}-\frac{k-1}{n-1}.
    \end{align*}

    For the case $x[n]=0$, it holds that   $j\notin \hat{e}$.   It follows from \reff{equ:jnotinE} that
    \begin{align*}
      \alpha_n(K_n^-) =  \sum_{e\notin E(j)}\mathcal{L}(e) x^k+ \sum_{e\in E(j)}\mathcal{L}(e) x^k \ge \binom{n-2}{k-2}.
    \end{align*}
    where the last inequality follows from the arithmetic-geometric mean inequality that $\mathcal{L}(e) x^k\ge0$ 
    for all $e\notin E(j)$. 
    In fact, the lower bound can be achieved by   set   $\tilde{x}$ as
    \begin{align*}
        \tilde{x}[i]=\left\{ \begin{array}{cl}
          (\frac{1}{n-1})^{\frac1k},&\quad \text{if} \quad i=1,\cdots,n-1;\\
          0,                         &\quad \text{if} \quad i=n.
        \end{array}
        \right.
    \end{align*}
    Hence,  $ \alpha_n(K_n^-)=\binom{n-2}{k-2}.$
 
%
%
    Hencefore,
    \begin{align*}
      \alpha_1(K_n^-) < \binom{n-2}{k-2} =\alpha_n(K_n^-).
    \end{align*}
    This complete the proof of \reff{equ:j1lj2}.
  \eproof



As discussed above, those vertices of the smallest degree are highly possible to help attain the analytic connectivity of a uniform hypergraph. A conjecture comes as follows.

\noindent\textbf{Conjecture 3.1}\quad Let $G=([n],E)$ be a $k$-graph. $\alpha(G)=\alpha_j(G)$ for some $j\in [n]$ of the smallest degree.


\section{A feasible trust region algorithm}
\label{section:FTRalg}

In this section, we propose the feasible trust region method ({\tt FTR}) for solving \reff{Def:anacon}.
Noting that the projection to the $k$-norm sphere and nonnegative space are easy.
Hence, we manage to  project the iterate points to the feasible set, while
maintaining the convergence.


The  problem \reff{Def:anaconj} can be rewritten  as follows
 \begin{align}\label{equ:maxforg}
   \nn \alpha_j=\min_{x\in\mbR^n} \quad&  \frac1{k}\mathcal{L} x^k,\\
  \nn \text{s.t.} \quad &\frac1k\left(\sum_{i=1}^nx[i]^k -1\right)=0, \\
                    & x\ge0,\quad x[j]=0, 
   \end{align}
which is equivalent to
  \begin{align}\label{equ:maxf}
   \nn \min_{x\in\mbR^{n-1}} \quad&  f(x)=\frac1{k}\tilde{\mathcal{L}} x^k,\\
  \nn \text{s.t.} \quad &c(x):=\frac1k\left(\sum_{i=1}^{n-1}x[i]^k -1\right)=0, \\
                    & x\ge0, 
   \end{align}
   where $\tilde{\mathcal{L}}\in T_{k,n-1}$ is the subtensor of $\mathcal{L}$ indexed by $[n]\setminus \{j\}$.


Before describing the details of {\tt FTR} algorithm, the following   functions are given.
The Lagrangian function of \reff{equ:maxf} is
\begin{align}\label{equ:L(x)}
 L(x,\lambda) = f(x) -  \lambda c(x),
\end{align}
and its gradient vector and Hessian matrix are
\begin{align}
\label{equ:g(x)}
  g(x)=\nabla_xL(x,\lambda)& = \nabla f(x) -  \lambda \nabla c(x),\\
  \label{equ:W(x)}
  W(x)=\nabla^2_{xx}L(x,\lambda)&=\nabla^2 f(x)-\lambda \nabla^2 c(x),
\end{align}
where $\nabla f(x)=\mathcal{\tilde{L}}x^{k-1}$, $\nabla^2 f(x)=(k-1)\mathcal{\tilde{L}}x^{k-2}$,
$\nabla c(x)=x^{[k-1]}$, $\nabla^2 c(x)=(k-1) \text{diag}(x^{[k-2]})$.
Here, $\mathcal{\tilde{L}}x^{k-1}\in\mbR^{n-1}$ is a vector with the $i$-th element being
\[(\mathcal{\tilde{L}}x^{k-1})[i] = \sum_{i_2,\cdots,i_k =1}^{n-1} \mathcal{\tilde{L}}_{i,i_2,\cdots,i_k}x_{i_2}\cdots,x_{i_k},\]
and $\mathcal{\tilde{L}}x^{k-2}\in\mbR^{(n-1)\times (n-1)}$ with the $(i, j)$-th element denoted as
\[(\mathcal{\tilde{L}}x^{k-2})[i, j] =
\sum_{i_3,\cdots,i_k =1}^{n-1} \mathcal{\tilde{L}}_{i, j, i_3,\cdots,i_k}x_{i_3}\cdots,x_{i_k}.\]
The  function vector $\mathcal{\tilde{L}}x^{k-1}$ is the subvetcor of $\mathcal{L}x^{k-1}$,
indexed  by  $[n]\setminus \{j\}$, and
the matrix   $\mathcal{\tilde{L}}x^{k-2}$ is $[n]\setminus \{j\}$ submatrix of $\mathcal{L}x^{k-1}$.

\subsection{The feasible trust region algorithm}

Given the current  point $x_t$,  the  trust region subproblem of \reff{equ:maxf} can be reformulated as follows,
\begin{align}\label{equ:orgmaxsub}
 \nn \min_{d\in\mbR^{n-1}} \quad &  m_t(d) = f_t + g_t^{\mathrm{T}}d + \frac12 d^{\mathrm{T}}W_t\,d,\\
  \nn \text{s.t.}\quad & c(x_t) + \nabla c(x_t)^{\mathrm{T}}d = 0,\\
   \nn      & \|d\| \le \Delta_t,\\
            & x_t +d\ge 0.
\end{align}
where $f_t=f(x_t)$,  $g_t=\nabla_xL(x_t,\lambda_t)$, $W_t=\nabla^2_{xx}L(x_t,\lambda_t)$,
$\Delta_t$ is the trust region radius updated in \reff{equ:Delta}.

In order to facilitate the computation of \reff{equ:orgmaxsub}, we utilize  the following strategies.
Firstly, we adopt the $\infty$-norm in \reff{equ:orgmaxsub},
and hence all the constrains will be linear.
Secondly, at each iteration, the feasibility of $x_t$ implies that $c(x_t)=0$, which ensures the feasibility of the resulting trust region subproblem.
Consequently, each subproblem is formulated as
\begin{align}\label{equ:maxsub} 
 \nn \min_{d\in\mbR^{n-1}} \quad &  m_t(d) = f_t + g_t^{\mathrm{T}}d + \frac12 d^{\mathrm{T}}W_t\,d,\\
  \nn \text{s.t.}\quad & \nabla c(x_t)^{\mathrm{T}}d = 0,\\
   \nn      & \|d\|_{\infty} \le \Delta_t,\\
            & x_t +d\ge 0.
\end{align}

Specifically,  at each iteration, if the trial step $d_t$ is accepted, the iterate $x_t+d_t$ is
projected to be feasible by setting $x_{t+1}=P(x_t+d_t)$, where
\begin{equation} P(x) = \frac{x}{\|x\|_k} \label{rx} \end{equation}
is a projection operator to the $k$-norm sphere and $\|x\|_k = (\sum_{i=1}^nx_i^k)^{1/k}$ is the $k$-norm of $x$. Set
\begin{align}\label{equ:lambdak}
\lambda_t=\nabla f(x_t)^{\mathrm{T}}x_t=\mathcal{A} x_t^m.
\end{align}  which is actually the Lagrange multiplier as will be clarified in \reff{equ:lmd}.

%


The following definitions are commonly used in trust region methods.
Denote the ratio of actual decrease and predicted decrease as
\begin{equation}\label{equ:rho}
  \rho_t = \frac{f(x_t) - f\left(P(x_t+d_t)\right)}{ m_t(0) - m_t(d_t)}.
\end{equation}
This is an important value for evaluating the error
between $m_t(d)$ and $f(x)$ at $x_t$. If $\rho_t$ is large enough,
we are confident to increase the trust region radius $\Delta_t$;
but if  $\rho_t$ is less than a threshold, we have to decrease the radius.
Specifically, $\Delta_{t+1}$ is updated as follows
\begin{equation}\label{equ:Delta}
  \Delta_{t+1} = \left\{ \begin{array}{ll}
 \frac12 \Delta_t,&\ \text{if}\ \rho_t \le \sigma_1; \\
                      \min\left(\Delta_{\max},2 \Delta_t\right),&\ \text{if}\ \rho_t>\sigma_2;\\
                      \Delta_t,&\ \text{else},
                    \end{array}
             \right.
\end{equation}
 where $\sigma_1, \sigma_2$ are constants with $0 < \sigma_1<\sigma_2 $ and $\sigma_1 < 1$.
 We only update $x_{t}$ in the next iteration  when $\rho_t$ is greater than or equal to some threshold,
\begin{equation}\label{equ:accept}
  x_{t+1} = \left\{ \begin{array}{ll}
 P(x_t + d_t),&\ \text{if}\ \rho_t \ge \sigma_0; \\
                      x_t,&\ \text{else},
             \end{array}
             \right.
\end{equation}
where $\sigma_0\in (0,\sigma_1)$ is a constant.
It should be noted that when updated,
$x_{t+1}$ is defined as the projection $P(x_t + d_t)$ instead of $x_t + d_t$.

The detailed descriptions of the {\tt FTR} method for computing the  analytic connectivity
\reff{Def:anacon} of symmetric tensors is as follows.
The algorithm includes two steps: the outer step and the inner step.
In the outer step, given an  index $j$,  let $x[j]=0$, and compute $\alpha(G)=\min_j \alpha_j(G)$.
In the inner step, the problem   \reff{equ:maxsub} is solved
by the feasible trust region algorithm to compute $\alpha_j(G)$.

\begin{algorithm}[H]\label{alg:FTR}
\caption{The feasible trust region method for the problem \reff{Def:anacon}}
\begin{description}
 \item[Step 0.] Given an initial point $x_0$, set the parameters $\sigma_0, \sigma_1, \sigma_2, \epsilon$, $\Delta_0$, $\Delta_{\max}$. Let $j=1$, $iter$ =0.
  \item[Step 1.] For  $ j=1,\cdots n$, do
  \begin{itemize}
 \item[s0.]     $\lambda_{0}=\mathcal{A} {x_{0}^m}$ and $t:=0$.
 \item[s1.] Solve the quadratic problem \reff{equ:maxsub} to determine $d_t$.
  \item[s2.]  If $\|d_t\| \le \epsilon$, stop and output $(\alpha_j(G)= \lambda_t$, $x^j =x_{t})$.
            Let $iter = iter+t$, and go to {\textbf{Step 1}}.
 \item[s3.] Calculate  $\rho_t$ by  \reff{equ:rho}.
 \item[s4.] Update the trust region radius $\Delta_t$ by \reff{equ:Delta}.
 \item[s5.] If $\rho_t \ge \sigma_0$, set $x_{t+1} = P(x_t+d_t)$ and $\lambda_{t+1}=\mathcal{A} {x_{t+1}^k}$;
 else $x_{t+1} = x_t$ and $\lambda_{t+1}=\lambda_t$. Set {$t$ $ := t+1$} and go to s0.
  \end{itemize}
\item[Step 2.] Let $j^* = \arg\min_{j=1}^n \alpha_j(G)$. Output $(\alpha_{j^*}(G), x^{j^*})$ and $iter$.
 \end{description}
\end{algorithm}

It is worth pointing out that if the involved uniform hypergraph has some special structure, such as those discussed in Section 3, then the computation in Algorithm \ref{alg:FTR} can be significantly reduced since the number of the outer loop can be cut down by merely considering those $j$ of the minimum degree.


\section{Convergence analysis}\label{section:convergence}

The first-order and the second-order optimality conditions of \reff{equ:maxf} are stated, and the global convergence of Algorithm 1 is established in this section.

\subsection{Optimality conditions}
For any local minimizer $x^*$ of \reff{equ:maxf}, the fact $\nabla c(x^*)=(x^*)^{[k-1]}$ implies 
that the set $\{\nabla c(x^*)\}\cup \{ e_i: i\in\mathcal{I}(x^*)\}$ is linearly independent,
where $e_i\in\mathbb{R}^n$ is the identity vector with the $i$-th element being one while the other elements are zero, 
and $\mathcal{I}(x^*)$ is the active set of $x^*$.
 Thus, the linear independence constraint qualification (LICQ) holds automatically. This observation immediately leads to the following first-order and second-order necessary conditions for \reff{equ:maxf} by invoking Theorems 12.1 and 12.5 in \cite{Nocedal2006}.
\begin{lemma}\label{lemma:1KKT}(First-order necessary conditions)
  Suppose that $x^*$ is a local solution of \reff{equ:maxf}.
  Then there is a Lagrange multiplier $\lambda^*$ such that
  \begin{equation}\label{equ:KKT}
    \min(x^*, g^*)=0,\quad c(x^*)=0,
  \end{equation}
  where $g^*=\nabla_x L(x^*,\lambda^*)=\nabla f(x)^*-\lambda^*\nabla c(x^*)$.
  Further, we have
  \begin{align}\label{equ:lmd}
    \lambda^*= (\nabla f(x)^*)^Tx^*.
  \end{align}
\end{lemma}


\begin{lemma}(Second-order necessary condition)
  Suppose that $x^*$ is a local solution of \reff{equ:maxsub}.
  Let $\lambda^*$ be the Lagrange multiplier satisfying \reff{equ:KKT}.
  Then
  \begin{equation}
    \label{equ:2KKT}
    d^T W^*d\ge0,\quad \forall\, d\in\,\mathcal{C}(x^*,\lambda^*),
  \end{equation}
  where  \begin{align}\label{equ:mathcalC}
    \mathcal{C}(x^*,\lambda^*)=\{d\mid\nabla c(x^*)^Td=0;\
    d[i]=0, \forall i\in\mathcal{I}(x^*)\ \text{with}\ g^*[i]>0;\
    d[i]\ge0,\forall  i\in\mathcal{I}(x^*)\ \text{with}\ g^*[i]=0\},
\end{align} and $W^*=\nabla_{xx}^2L(x^*,\lambda^*)$.
\end{lemma}

\subsection{Global convergence}
In this subsection, we establish the global convergence  of
 the inner problem of Algorithm \ref{alg:FTR}; i.e., using
feasible trust region algorithm to solve the  problem \reff{equ:maxf}.
 We shall employ the techniques in traditional trust region methods to derive the results.
However, there are two key difficulties.
 Firstly,   $x_{t+1}$ is updated by $P(x_t+d_t)$ instead of $x_t+d_t$ in order to keep the feasibility.
We should estimate the error between $f(P(x_t+d))-f(x_t)$ with its second order approximation,
instead of $f(x_t+d)-f(x_t)$. Secondly, $\infty$-norm is applied, hence the outline of proof is
different from Euclidean-norm cases.

 To simplify our analysis, define
 \[h(x)=f(P(x)).\]
Then the gradient and the Hessian of $h(x)$ are
 \begin{align}
  \nabla h(x)&=\nabla P(x)\nabla f(P(x)),\label{equ:h1}
    \end{align}
     \begin{align*}
 \nabla^2 h(x) &=\frac{\nabla P(x)\nabla^2 f(P(x))}{\|x\|_k}
           - \frac{\nabla c(x) \nabla f(P(x))^{\mathrm{T}}}{\|x\|_k^{k+1}}
           +\frac{(k+1)\nabla x^{\mathrm{T}}f(P(x))\nabla c(x) \nabla c(x)^{\mathrm{T}}}{\|x\|_k^{2k+1}}\\ 
    & -\frac{x^{\mathrm{T}}\nabla f(P(x))\nabla^2 c(x)+ \nabla P(x)\nabla^2 f(P(x))x\nabla c(x)^{\mathrm{T}}
    + \nabla f(P(x))\nabla c(x)^{\mathrm{T}}}{\|x\|_k^{k+1}},
 \end{align*}
 where $\nabla P(x)=\left(\frac{I}{\|x\|_k} - \frac{\nabla c(x)x^{\mathrm{T}}}{\|x\|_k^{k+1}}\right)$.
A key property is that when $\|x_t\|_k=1$ and $\nabla c(x_t)^{\T}d=0$, we have
\begin{equation}\label{equ:nablah}
  \nabla h(x_t)^{\T}d =\nabla f(x_t)^{\T}d= g(x_t)^{\T}d
\end{equation}
and
\begin{equation}\label{equ:nabla2h}
  d^{\T}\nabla^2 h(x_t)d=d^{\T}\nabla^2 f(x_t)d-\lmd d^{\T}\nabla^2 c(x_t)d=d^{\T}W(x_t)x.
\end{equation}
That is, the feasible direction $d$ satisfying $\nabla c(x)^Td=0$,
the second order approximations of $h(x)$ and $L(x, \lmd)$ are the same. Several
technical lemmas are presented for the convergence analysis.


\begin{lemma}\label{lem:gB}

\begin{itemize}
  \item[(i)]  Let $g(x)$ and $W(x)$ are defined in \reff{equ:g(x)} and \reff{equ:W(x)}, respectively.
  When $\lambda$ is fixed, for all $x\ge0$, $y\ge0$ satisfying $\|x\|_k=1$ and $\|y\|_k=1$,  we have
 \begin{align}
\|W(x)\|&\le M,                    
    \label{equ:B<M}\\
    \|g(x)-g(y)\|&\le L_0\|x-y\|,&
    \label{equ:gM0}\\
    \|W(x)-W(y)\|&\le L_1\|x-y\|,&
    \label{equ:Bx-By<M}
 \end{align}
   where $M$, $L_0$ and $L_1$  are positive constants.
      \item[(ii)] Suppose $\|x\|_k\ge \eta_1$, $\|y\|_k\ge \eta_2$,
      where $\eta_1$ and $\eta_2$ are positive constants.  We have
     \begin{equation} \label{equ:L2}
     \| \nabla^2 h(x) - \nabla^2 h(y)\| \le L_2\|x-y\|,
     \end{equation}
  where $L_2$  is a positive constant.
     \end{itemize}

\end{lemma}
\textbf{Proof.} They are obvious since  $g(x)$, $W(x)$ and $\nabla^2 h(x)$
are {smooth and bounded} on the closed sets.
\quad \eproof

\begin{lemma}\label{lem:md-fd}
Suppose $x_t$ is feasible solution of model \reff{equ:maxf},
and  $d_t$ is feasible solution of model \reff{equ:maxsub}.
 For the error between the models {$m_t(d_t)$ }and $h(x_t+d_t)$, we have
  \begin{equation}\label{equ:md-fd}
    \left|m_t(d_t) - h(x_t+d_t) \right| \le \beta \|d_t\|^3,
  \end{equation}
  where  $\beta$ is some positive constant. 
\end{lemma}
\textbf{Proof.}\ By the mean value theorem for integration, we have
\begin{align*}
  h(x_t+d_t)& =  h(x_t) + \nabla h(x_t)^{\mathrm{T}}d_t + \frac12d_t^{\mathrm{T}}\nabla^2 h(x_t+\theta_t d_t)d_t
\end{align*}
for some $\theta_t \in(0,1)$.
It follows from $h(x_t)=f(x_t)$, \reff{equ:nablah}  and \reff{equ:nabla2h} that
\begin{align*}
 \left| m_t(d_t) - h(x_t+d_t)\right| &= \left|\frac12 d_t^{\mathrm{T}}W_td_t -  \frac12d_t^{\mathrm{T}}\nabla^2 h(x_t+\theta_t d_t)d_t\right|\\
                                 &= \left|\frac12 d_t^{\mathrm{T}}\nabla^2 h(x_t)d_t -  \frac12d_t^{\mathrm{T}}\nabla^2 h(x_t+\theta_t d_t)d_t\right|\\
                                     &\le \frac12L_2 \|d_t\|^3.
 \end{align*}
To show the above inequality by Lemma \ref{lem:gB} ($ii$), we still need to prove
 $\|x_t\|_k$ and $\|x_t+\theta d_t\|_k$ are positive.
The feasible point  $x_t$ satisfies $\|x_t\|_k=1$.
 As  two norms are equivalent, i.e., for $x\in\mathbb{R}^n$ if $r_1>r_2>0$, then
  \[\|x\|_{r_1}\le\|x\|_{r_2}\le n^{\frac1{r_2}-\frac1{r_1}}\|x\|_{r_1}. \]
 Hence, it follows from  $\|x_t\|_k=1$ that for $k\ge3$, $ \| \nabla c(x_t)\|=\|x_t^{[k-1]}\|=\|x_t\|_{2k-2}^{k-1}\le 1$.
Furthermore,  it follows from  $\nabla c(x_t)^{\mathrm{T}}d_t=0$,
  $\nabla c(x_t)=x_t^{[k-1]}$ and $\nabla c(x_t)^{\mathrm{T}}x_t=\|x_t\|_k^k=1$ that
 $ \nabla c(x_t)^{\mathrm{T}} (x_t+\theta_t d_t)=1$.
 Therefore,
 \begin{equation*}
 \| \nabla c(x_t)\|\cdot \|x_t+\theta_t d_t\| \ge 1.
  \end{equation*}
As a result,  both $x_t$ and $x_t+\theta_td_t$ are lower bounded.
\quad\eproof

%
%

\begin{lemma}\label{lem:infg=0}
  Consider the sequence  $\{x_t\}$ generated by Algorithm \ref{alg:FTR}. Then sequence  $\{f(x_t)\}$
    of  the objective value is  nondecreasing. Furthermore,  at least one of the cluster points of $\{x_t\}$ is a KKT
  points of the problem \reff{equ:maxf}, i.e.,
 \begin{equation}\label{equ:convergence}
   \liminf_{t\rightarrow \infty}\|\min(x_t, \nabla f(x_t)-\lambda_t\nabla c(x_t))\| =0.
 \end{equation}
 \end{lemma}

\textbf{Proof.}
Suppose   the theorem is false,  we assume that
\begin{equation}\label{equ:Delgoto0}
  \lim_{k\rightarrow\infty}\Delta_t =0.
\end{equation}
If \reff{equ:Delgoto0} fails, there exists a const $\delta>0$, such that for infinite many $t$,
it holds that
\begin{equation}\label{equ:contrdict}
  \Delta_t\ge\delta \quad \text{and}\quad \rho_t\ge\sigma_1.
\end{equation}
Denote the set of $k$ satisfying  \reff{equ:contrdict} as $K_0$.
Without loss of generality, suppose
\begin{equation*}
  \lim_{t\in K_0,  t\rightarrow \infty} x_t=\bar{x}.
\end{equation*}
 According to our assumption, $\bar{x}$ is not a stationary point of \reff{equ:maxf},
 hence $d=0$ is not the optimal solution of the following system
 \begin{align}\label{equ:contrsubpro}
 \nn \min_{d\in\mbR^{n-1}} \quad &  \bar{m}(d)= f(\bar{x}) + g(\bar{x})^{\mathrm{T}}d + \frac12 d^{\mathrm{T}}W(\bar{x})\,d,\\
  \nn \text{s.t.}\quad & \nabla c(\bar{x})^{\mathrm{T}}d = 0,\\
   \nn      & \|d\|_{\infty} \le \delta,\\
            & \bar{x} +d\ge 0.
\end{align}
Denote $\bar{d}$ as its solution, then
\begin{equation*}
    \gamma =\bar{m}(0)-\bar{m}(\bar{d})= -g(\bar{x})^{\mathrm{T}}\bar{d} - \frac12 \bar{d}^{\mathrm{T}}W(\bar{x})\,\bar{d}>0.
\end{equation*}
It follows from Lemma \ref{lem:subprobcont} that
$$m_t(0)-m_t(d_t)\ge \frac{1}{2}(\bar{m}(0)-\bar{m}(\bar{d}))\ge\frac12\gamma$$
for all $t\in K_0$ large enough. As a result,
$f(x_t)-f(x_{t+1})\ge\frac12\sigma_1\gamma >0$
for all large enough $t\in K_0$.
This contradicts to $\lim_{t\rightarrow\infty}f(x_t)=f(\bar{x})$.
The contradiction indicates that \reff{equ:Delgoto0} holds.

If \reff{equ:Delgoto0} holds, there exists a subsequence such that
\begin{equation}\label{equ:contr2}
    \rho_t\le\sigma_1,\quad\forall\, t\in K_1.
\end{equation}
Without loss of generality, suppose
\begin{equation}\label{equ:xgotohatx}
  \lim_{t\in K_1, t\rightarrow\infty} x_t =\hat{x},
\end{equation}
 According to our assumption, $\hat{x}$ is not a stationary point of \reff{equ:maxf},
 hence $d=0$ is not the optimal solution of the following system
 \begin{align}\label{equ:contrsubpro2}
 \nn \min_{d\in\mbR^{n-1}} \quad &  \hat{m}(d)= f(\hat{x}) +g(\hat{x})^{\mathrm{T}}d + \frac12 d^{\mathrm{T}}W(\hat{x})\,d\\
  \nn \text{s.t.}\quad & \nabla c(\hat{x})^{\mathrm{T}}d = 0,\\
   \nn      & \|d\|_{\infty} \le 1,\\
            & \hat{x} +d\ge 0.
\end{align}
Denote $\hat{d}$ as its solution, then
\begin{equation*}
    \hat{\gamma} =\hat{m}(0)-\hat{m}(\hat{d})= -g(\hat{x})^{\mathrm{T}}\hat{d} - \frac12 \hat{d}^{\mathrm{T}}W(\hat{x})\,\hat{d}>0.
\end{equation*}
As $\hat{d}_t=\Delta_t \hat{d}$ is the solution of \reff{equ:contrsubpro2} with the
trust region radius replaced by $\Delta_t$. Then $\hat{m}(0)- \hat{m}(\hat{d}_t)\ge \frac12\Delta_t\hat{\gamma}$.
It follows from Lemma \ref{lem:subprobcont} that
\begin{equation}\label{equ:contrdecre}
  m_t(0) - m_t(d_t) \ge\frac12( \hat{m}(0)- \hat{m}(\hat{d}_t))\ge\frac14\Delta_t\hat{\gamma}
\end{equation}
for all $t\in K_1$ large enough, where the last inequality comes from $\Delta_t\rightarrow 0$.
Further,
\begin{align}\label{equ:rhoged}
 \nn   \rho_t&\ge1-|1-\rho_t|\\
  \nn          &=1-\frac{|m_t(0)-m_t(d_t)+h(x_t+d_t)-h(x_t)|}{|m_t(0)-m_t(d_t)|}\\
  \nn         &=1-\frac{|h(x_t+d_t)-m_t(d_t)|}{|m_t(0)-m_t(d_t)|}\\
                &\ge1-\frac{\beta \|d_t\|^3}{|m_t(0)-m_t(d_t)|}.
\end{align}
This, together with \reff{equ:contrdecre}, derives $\lim_{t\in K_1, t\rightarrow \infty}\rho_t=1$,
which contradicts  with \reff{equ:contr2}. This completes the proof.
\quad \eproof

\begin{lemma}\label{lem:subprobcont}

The optimal value of \reff{equ:contrsubpro} is continuous
for all feasible points $\bar{x}$ of  \reff{equ:maxf}. Namely,
given two points $x_{t_1}$ and $x_{t_2}$ satisfying
$\|x_{t_i}-\bar{x}\|\le\epsilon_0,\ i=1,2$
with $x_{t_i}\ge0, \ \|x_{t_i}\|_k=1,\ i=1,2$,
their optimal solution for  \reff{equ:contrsubpro} are
$d_{t_1}$ and $d_{t_2}$, respectively. Then, for all $\epsilon>0$ small enough, it holds that
\begin{align}\label{equ:lemma4.4}
|g(x_{t_1})^\T d_{t_1}+\frac12d_{t_1}^\T W(x_{t_1})d_{t_1} -g(x_{t_2})^\T d_{t_2}
-\frac12d_{t_2}^\T W(x_{t_2})d_{t_2}|\le \epsilon.
\end{align}

\end{lemma}
\textbf{Proof.}\quad
As $\bar{x}$ satisfies $\sum_i \bar{x}[i]^k=1$, there exists at least an index $p$ such that $\bar{x}[p]>0$.
For two points $x_{t_1}$ and $x_{t_2}$ near $\bar{x}$, there exists a positive value $\epsilon_1$ such that
$$\|\nabla c(x_{t_2})^{\mathrm{T}}d_{t_1}\|=
\|\nabla c(x_{t_2})^{\mathrm{T}}d_{t_1}-\nabla c(x_{t_1})^{\mathrm{T}}d_{t_1}\|\le \|d_{t_1}\|\|\nabla c(x_{t_2})-\nabla c(x_{t_1})\|\le\epsilon_1,$$
where the last inequality follows from that $d_{t_1}$ is bounded, and $\nabla c(x)$ is continuous.
If $d_{t_1}[p]<\delta$ and $\nabla c(x_{t_2})^{\mathrm{T}}d_{t_1} >0$
or $d_{t_1}[p]>-\delta$ and $\nabla c(x_{t_2})^{\mathrm{T}}d_{t_1} <0$,
then $\tilde{d}_{t_2}=d_{t_1}-\frac{\nabla c(x_{t_2})^{\mathrm{T}}d_{t_1}}{\nabla c(x_{t_2})^{\mathrm{T}}e_p}e_{p}$ is a feasible solution for
  \begin{align}\label{equ:contasubprok2}
 \nn \chi_{t_2}= \min_{d\in\mbR^n} \quad & g(x_{t_2})^{\mathrm{T}}d + \frac12 d^{\mathrm{T}}W(x_{t_2})\,d,\\
  \nn \text{s.\,t.}\quad & \nabla c(x_{t_2})^{\mathrm{T}}d = 0,\\
   \nn      & \|d\|_{\infty} \le \delta,\\
            & x_{t_2} +d\ge 0.
\end{align}
Otherwise, suppose that $d_{t_1}[p]=\delta (-\delta)$, from $\nabla c(x_{t_1})^{\mathrm{T}}d_{t_1}=0$
that there exists some positive index $q$ such that $\bar{x}[q]>0$ and $d_{t_1}[q]<(>)0$,
hence $\tilde{d}_{t_2}=d_{t_1}-\frac{\nabla c(x_{t_2})^{\mathrm{T}}d_{t_1}}{\nabla c(x_{t_2})^{\mathrm{T}}e_q}e_{q}$ is a feasible solution for the above problem.
Therefore, from the fact that the objective function of \reff{equ:contasubprok2} is continuous
and that $\tilde{d}_{t_2}$ is only a feasible solution, we have
$$\chi_{t_2}\le g(x_{t_2})^{\mathrm{T}}\tilde{d}_{t_2} + \frac12 \tilde{d}_{t_2}^{\mathrm{T}}W(x_{t_2})\,\tilde{d}_{t_2}\le\chi_{t_1}+\epsilon.$$
On the other hand, we can show $\chi_{t_2}+\epsilon\ge\chi_{t_1}$.
Therefore, \reff{equ:lemma4.4} holds true.
\quad \eproof

%

\begin{theorem}\label{thm:converge2orderpoint}
Suppose that  the iterates $\{x_t\}$ generated by Algorithm \ref{alg:FTR} converge to  $x^{*}$.
Then the second-order necessary conditions \reff{equ:2KKT} holds.
\end{theorem}
\textbf{Proof.}\quad
We show this theorem  by contradiction.
Suppose that there exists a negative eigenvalue $-\eta_0$ satisfying
\begin{equation}\label{equ:vBv>0}
  v^{\text{T}}W^*v=-\eta_0 <0,\quad\text{where}\  v\in \mathcal{C}(x^*,\lambda^*),\quad \|v\|_2=1.
\end{equation}
It follows from the definition of \reff{equ:mathcalC} that $v$ is a feasible solution of
\reff{equ:maxsub} with $x_t$ replaced by $x^*$, and $\Delta_t$ replaced by 1. 
For all $i\in  \mathcal{I}(x^*)$, either $g^*[i]=0$ or $v[i]=0$,
and for all $i\notin \mathcal{I}(x^*)$, $g^*[i]=0$, hence
 \begin{equation}\label{equ:gv=0}
  (g*)^\T v=0.
\end{equation}
When $x_t$ is close enough to $x^*$, it follows from the proof of Lemma \ref{lem:subprobcont}
and  that
 $\hat{d}_t =  \Delta_t v+d_t^{\epsilon}$ is a feasible point for the problem  \reff{equ:maxsub},
where  $\|d_t^{\epsilon}\|$ is small enough to  be bounded by $\|x_t-x^*\|$.
Furthermore, it follows from \reff{equ:gv=0} that $g_t^{\text{T}}\hat{d}_t$ is small,
  $v^{\text{T}}W^*v=-\eta_0 <0$. Hence $\hat{d}_t$
is an decrease direction for the problem  \reff{equ:maxsub}. Therefore,
\begin{align}
 \nn m_t(0)-m_t(\hat{d}_t)  &= -g_t^{\text{T}}\hat{d}_t - \frac12 d^{\text{T}}W_t\hat{d}_t\\
 \nn &= -\Delta_t  g_t^{\text{T}}v -\frac12\Delta_t^2v^{\text{T}}W_tv + o(1)\\
\nn &\ge- \frac14\Delta_t^2v^{\text{T}}W_tv.
\end{align}
 Since  $|v^{\text{T}}W_tv-v^{\text{T}}W^*v|\leq \|W_t-W^*\|\|v\|^2$, $\|v\|=1$,
then $v^{\text{T}}W_tv\le-\frac12 \eta_0$ and
\begin{equation}\label{equ:md-m0>lmd}
   m_t(0)-m_t(d_t)\ge   m_t(0)-m_t(\hat{d}_t)\ge-\frac14\Delta_t^2v^{\text{T}}W_tv\ge \frac18\Delta_t^2\eta_0.
\end{equation}
It follows from \reff{equ:rhoged} that $\rho_t\rightarrow 1$.
Therefore, there exists $K_2$ large enough such that
\begin{equation}
  f(x_t)-f(x_{t+1})  \ge \sigma_1(m_t(0)-m_t(d_t)) \ge {\frac18\Delta_t^2\sigma_1\eta_0,} \quad \forall\, k\ge K_2,
\end{equation}
which derives that $\Delta_t\rightarrow0$. This contradicts with $\rho_t\rightarrow 1$.
Thus, \reff{equ:vBv>0} is false.
 \quad\eproof


\section{Numerical experiments}
\label{section:numericalexperiments}

In this section, we present several  numerical results of computing the analytic connectivity. Our codes are implemented in MATLAB (R2014a).
All the experiments are preformed on a Dell desktop with Intel dual core i7-4770 CPU at
3.40 GHz and 8GB of memory running  Windows 7.
The parameters are set as
 \begin{equation*}
   \sigma_0 =0.25,\ \ \sigma_1 =0.5,\ \ \sigma_2 =0.75,\ \ \epsilon = 1.0^{-8},
   \ \ \Delta_0 = 2, \ \ \Delta_{\max} = 10.
 \end{equation*}
We execute the {\tt FTR}  algorithm 100 times with different initial points, and report the average results.
The  initial points are generated  by the following {\tt Matlab} commands
\begin{verbatim}
    for rd = 1:100;    randn('seed', rd);   x0 = randn(n-1,1);    end;
\end{verbatim}
which obey the Gaussian distribution.
Afterwards, $x_0$  is restricted to the feasible set of \reff{equ:maxf} by doing the projection $P(|x_0|)$.

{\tt FTR} is compared with an Sparse Nonlinear OPTimizer solver {\tt SNOPT} \cite{SNOPT}, which is called by the free trial software {\tt TOMLAB}
\footnote{http://tomopt.com/tomlab/}.
The exact gradient and the Hessian are provided for {\tt FTR} and {\tt SNOPT}, and
both the quadratic programming subproblems of {\tt FTR} and {\tt SNOPT} are computed by {\tt SQOPT}.
Furthermore, for small dimensional problems,
 we utilize the global optimization software {\tt GloptiPoly 3}   \cite{henrion2009gloptipoly}
\footnote{http://homepages.laas.fr/henrion/software/gloptipoly/}
 to solve \reff{equ:maxf}, which can help us to judge whether our solution is the global optimal solution.
{\tt GloptiPoly 3}  relaxes the polynomial problem into a hierarchy of semidefinite subproblems,
which are solved by {\tt SDPNAL+} \cite{yang2015sdpnal+}.

Noting that the main computation of {\tt FTR} includes calculating
 $\mathcal{L}x^k$, $\mathcal{L}x^{k-1}$ and $\mathcal{L}x^{k-2}$.
  To deal with this, we adopt the methods in Chang, Chen and Qi \cite{CCQ}
  to calculate  $\mathcal{L}x^k$, $\mathcal{L}x^{k-1}$,
   where they   store a uniform hypergraph by a compact matrix
$G_r\in\mathbb{R}^{m\times k}$, where $m$ is the number of edges, and $k$ is the number
of vertices in an edge; namely, the $i$-th edges of the hypergraph is the $i$-th row of $G_r$ as
\[G(i,:) = (v_{i_1}, \cdots, v_{i_k}).\]
The computational method for $\mathcal{L}x^{k-2}$ follows the same strategy. Thus, the computation cost for  $\mathcal{L}x^k$, $\mathcal{L}x^{k-1}$,
$\mathcal{L}x^{k-2}$ are $O(mk)$, $O(mk^2+mnk)$ and $O(mk^3 + mn^2k^2)$, respectively.
It should also be noted that the sparsity ratio of $\mathcal{L}x^{k-2}$ is
\begin{align*}
  nnz(\mathcal{L}x^{k-2})=O\left(\frac{mk^2}{n^2}\right).
\end{align*}
Thus our method enjoys fast computation when the sparsity property is utilized.

\subsection{Comparison of {\tt FTR} with {\tt SNOPT} and {\tt GloptiPoly 3} for small size hypergraphs}
\begin{figure}[h]
  \centering
  \includegraphics[width=.8in]{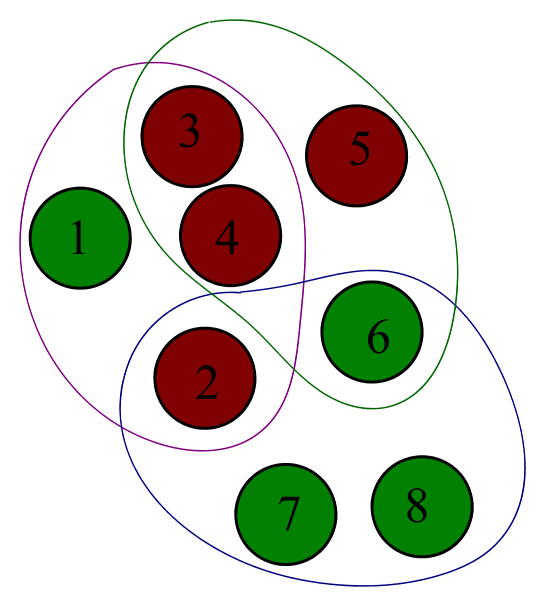}
  \includegraphics[width=1in]{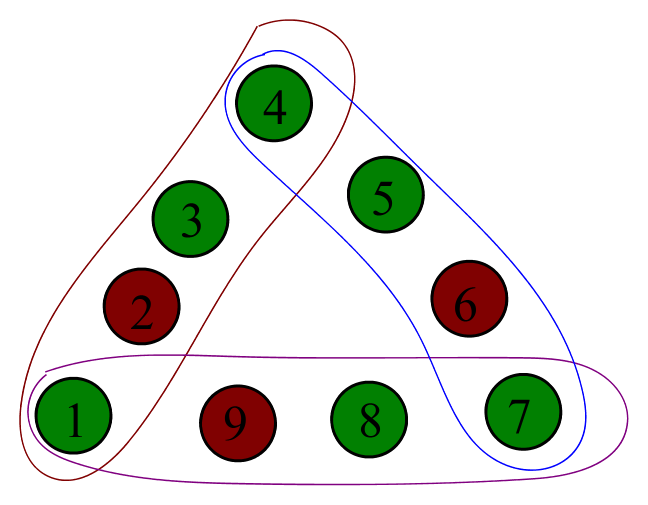}
  \includegraphics[width=.8in]{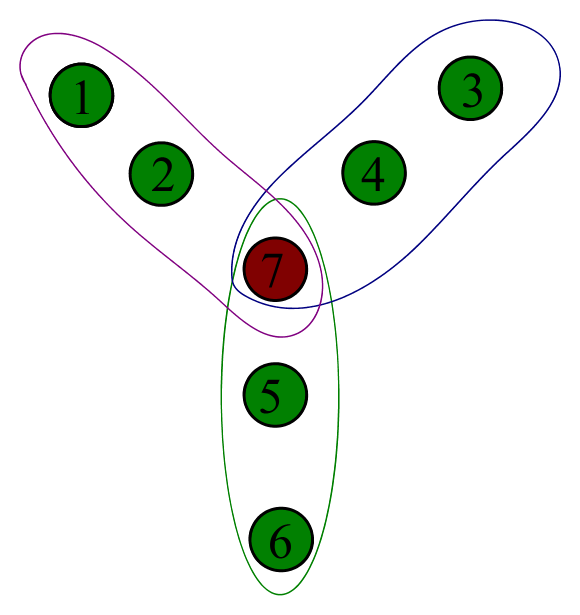}
   \includegraphics[width=.8in]{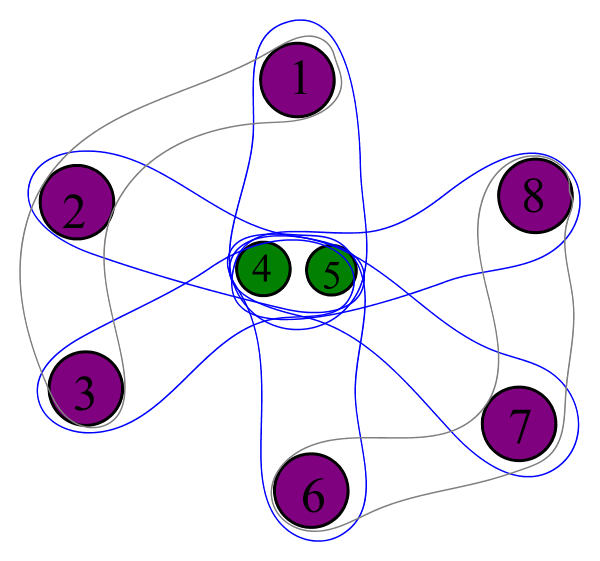}
    \includegraphics[width=1in]{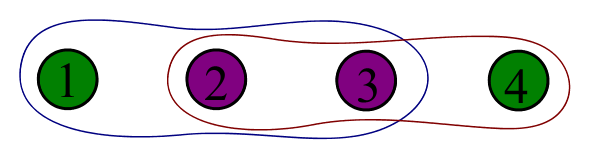}
     \includegraphics[width=1in]{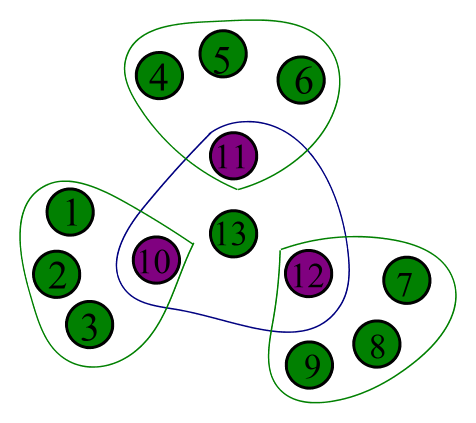}
  \caption{The uniform hypergraphs (a)--(f)}\label{fig:example}
\end{figure}

In this subsection, we show the numerical results of our {\tt FTR} algorithm, compared with {\tt SNOPT} and {\tt GloptiPoly 3}.
We will use the hypergraphs in Figure \ref{fig:example} which are found in \cite{CCQ, HQS, HQX, Qi} as the testing instances.

In Table \ref{HGanacontable1}, `$m$' is the number of edges of the hypergraph, `$n$' is the number of vertices,
$k$ is the number of vertices in an edge.
 `$\alpha$' means the analytic connectivity returned by {\tt FTR} and {\tt SNOPT},
`$\alpha^*$' stands for the analytic connectivity computed from the global optimization software {\tt GloptiPoly 3},
`ratio' means the ratio {\tt FTR} and {\tt SNOPT} get the same result with {\tt GloptiPoly 3},
and `iter' is the average number of iterations  of 100 runs with random initializations.
`time (s)' denotes the average CPU time of seconds consumed   by {\tt FTR} and {\tt SNOPT},
or the total CPU time of {\tt GloptiPoly 3}.

\begin{table}[h]
\centering
\caption{Comparisons of  {\tt FTR}  with {\tt SNOPT} and {\tt GloptiPoly 3}}
  \begin{tabular}{|c@{\hspace{1mm}}|@{\hspace{1mm}}c@{\hspace{2mm}}
  |@{\hspace{2mm}}c@{\hspace{1mm}}c@{\hspace{1mm}}r@{\hspace{1mm}}r@{\hspace{2mm}}
  |@{\hspace{2mm}}c@{\hspace{1mm}}c@{\hspace{1mm}}r@{\hspace{1mm}}r@{\hspace{2mm}}
  |@{\hspace{2mm}}c@{\hspace{1mm}}r@{\hspace{1mm}}|}
\hline
                         &          &    \multicolumn{4}{c}{{\tt SNOPT}}
                          &    \multicolumn{4}{@{\hspace{-0.2cm}}|c}{{\tt FTR}}
                          &  \multicolumn{2}{@{\hspace{-0.2cm}}|c|}{{\tt GloptiPoly 3}} \\ \hline
      Hypergraph           &($m, n, k$)&$\alpha$   & ratio  & iter & time (s)&$\alpha$   & ratio  & iter & time (s)  & $\alpha^*$& time (s)\\ \hline
   (a)&(3, 8, 4) &0.2516 &  100\% &323.61 &   0.3065   &0.2516 &  100\% & 75.42 &   0.0332   &0.2516 & 59.515          \\ \hline
  (b)&(3, 9, 4)  &0.2100 &  100\% &403.64 &   0.3619     &0.2100 &  100\% &  83.65 &   0.0365    &0.2100 &  110.14 \\ \hline 
  (c)& (3, 7, 3) &0.1607 &  100\% &142.15 &   0.1007      &0.1607 &  97\% &  48.15 &   0.0185      &0.1607 & 74.136 \\ \hline 
  (d)& (8, 8, 3) &0.4300 &  100\% &151.46 &   0.1216   &0.4300 &  100\% &67.03 &   0.0263   &0.4300 &  110.10 \\ \hline 
  (e)&(2, 4, 3)  &0.5344 &  100\% &42.28 &   0.0381  &0.5344 &  100\% & 25.28 &   0.0080  &0.5344 & 23.052 \\ \hline 
  (f)& (4, 13, 4)&0.0592 &  100\% &850.18 &   0.8496   &0.0592 &  97\% &131.77 &   0.0603 &0.0592 &   18.877 \\ \hline 
\end{tabular}
 \label{HGanacontable1}
\end{table}

Table \ref{HGanacontable1}  shows that  both {\tt SNOPT} and {\tt FTR} produce the
same results with {\tt GloptiPoly3} for almost 100\%.
This is in accord with Theorem \ref{thm:converge2orderpoint} that {\tt FTR} converges to
second order necessary points, which has a high possibility to converge to global optimal point.
Besides, the average iteration number that {\tt FTR} takes is relatively small comparing to that of {\tt SNOPT},
since {\tt FTR} has utilized the trust region technique. As the main computation costs in each iteration for both {\tt FTR} and {\tt SNOPT} are to solve the quadratic programming, this makes {\tt FTR} take less CPU time than {\tt SNOPT}, as one can see from Table \ref{HGanacontable1}. Additionally, it is known from Table \ref{HGanacontable1} that, among the above six hypergraph instances, the hypergraph (f) has the smallest analytic connectivity, while (d) and (e) have relatively large ones. This, to some extent, reflects the connectivity of the corresponding hypergraphs as can be seen from Figure \ref{fig:example}.

\subsection{Larger dimensional problems}

In this subsection, we are ready to compute relatively large dimensional problems by {\tt FTR}, and compare its performance with that of {\tt SNOPT}.
As {\tt GloptiPoly 3} will be too costly both in time and in space for large problems, we will not consider this algorithm here. Similar to the small dimensional cases, we also give 100 initial points, and show the overall and average results. We take the $2$-path $4$-graph as discussed in Proposition \ref{lemma:chaindecrease} and $K_n-$ in Proposition \ref{lem:kn-} for testing instances with different values of $n$. The computational results are shown in Tables \ref{HGanacontable10} and \ref{HGanacontablecomplete}, where `$\alpha$' is the analytic connectivity in question, and `ratio' stands for the percentage from 100 experiments to achieve that minimal value.
\vskip 2mm

\begin{table}[h]
\centering
\caption{Results for the $2$-path $4$-graphs with different $n$  by {\tt FTR} and {\tt SNOPT}}
  \begin{tabular}{|c|ccrr|ccrr|}
\hline
        &    \multicolumn{4}{c|}{{\tt FTR}}    &    \multicolumn{4}{c|}{{\tt SNOPT}}    \\ \hline
     $n$&$\alpha$  & ratio  & iter & time (s)&$\alpha$   & ratio  & iter & time (s) \\ \hline
        10 &1.21e-01  & 100\% &     11.67 &   0.0058 & 1.21e-01 & 100\% &     48.66 &   0.0408   \\ \hline
        50 &4.11e-03  &  92\%&     12.46 &   0.0095 & 4.11e-03 &  95\% &    186.16 &   0.1426      \\ \hline
        100&1.01e-03  &  82\% &     15.00 &   0.0233 & 1.01e-03 &  86\% &    268.11 &   0.2554      \\ \hline
        200&2.49e-04  &  98\% &     14.92 &   0.0872 & 2.49e-04 &  79\% &    534.81 &   1.3374     \\ \hline
        300&1.10e-04  &  95\% &     14.86 &   0.2274 & 1.10e-04 &  73\% &    816.92 &   4.6972     \\ \hline
        400&6.20e-05  &  96\% &     14.50 &   0.4935 & 6.20e-05 &  87\% &   1039.38 &  11.7781    \\ \hline
        500&3.96e-05  &  94\% &     14.71 &   0.9096 & 3.96e-05 &  89\% &   1329.87 &  26.4040     \\ \hline
\end{tabular}
 \label{HGanacontable10}
\end{table}

We can see from Table \ref{HGanacontable10} that both {\tt FTR} and {\tt SNOPT} produce the same optimal value for each of the above instances,
and the successful ratio is above 70\%, while {\tt FTR} is slightly better than {\tt SNOPT}. Comparing to those small size problems as computed in Subsection 6.1, 
large dimensional problems here are relatively hard to achieve the global optimum with local optimal algorithms such as {\tt FTR} and {\tt SNOPT}. For the iteration number, we find that  {\tt FTR} scales well for dimension as large as 500,  while {\tt SNOPT} takes far more iteration steps for larger dimensional problems. This leads to overwhelming superiority of {\tt FTR} in computation time comparing to {\tt SNOPT}, as one can see from Table \ref{HGanacontable10}. Besides, it is worth pointing out that the sparse ratio of the Hessian matrix for this problem
is about $O(\frac{1}{n})$, and both the quadratic subproblems of {\tt FTR} and {\tt SNOPT} have taken this advantage. Thus, the overall computation time is not long even when the iteration number as big as more than 1000. Additionally, we can see that as $n$ increases, $\alpha(G)$ is monotonically decreasing, which fits the result in Proposition \ref{lemma:chaindecrease}.

\vskip 2mm
The numerical results for $K_n^{-}$ with $k=3$ and different values of $n$ are shown in Table \ref{HGanacontablecomplete}, with the comparison on performances of {\tt FTR} and {\tt SNOPT}, and the upper bounds $\bar{\alpha}=n-2-\frac{2}{n-1}$ given in \reff{upperbound}. As already known from Proposition \ref{lem:kn-}, $\alpha(K_n^{-}) = \min_{j=1,\cdots,k}\alpha_j(K_n^{-})$. Combining with the inherited symmetric structure of $K_n^{-}$, we only need to compute $\alpha_1(K_n^{-})$.

\begin{table}[h]
\centering
\caption{Numerical results for $K_n^-$ with different $n$ by {\tt FTR} and {\tt SNOPT}}

  \begin{tabular}{|c|ccrr|ccrr|c|}
\hline
         &    \multicolumn{4}{c|}{{\tt FTR}}    &    \multicolumn{4}{c|}{{\tt SNOPT}}  & upper~bound   \\ \hline
           $n$&$\alpha$   & ratio  & iter & time (s)&$\alpha$   & ratio  & iter & time (s) & $\bar{\alpha}$  \\ \hline
      10&  7.7736  &100\% &    6.82  &   0.0031 &   7.7736  &100\% &  30.01  &   0.0262   &     7.7778 \\ \hline
      20& 17.8943  &100\% &    7.27 &   0.0072 &  17.8943  &100\% &   17.56 &   0.0226   &    17.8947 \\ \hline
      30& 27.9309  &100\% &    8.03 &   0.0242 &  27.9309  &100\% &   14.48 &   0.0455   &    27.9310\\ \hline
      40& 37.9487  &100\% &    8.67 &   0.0764 &  37.9487  &100\% &   13.29 &   0.1578   &    37.9487 \\ \hline
      50& 47.9592 &100\%  &    8.54 &   0.2082 &  47.9592 &100\%  &   14.72 &   0.5159   &    47.9592    \\ \hline
      60& 57.9661 &100\%  &    8.38 &   0.4900 &  57.9661 &100\%  &   15.18 &   1.8829   &    57.9661     \\ \hline
      70& 67.9710 &100\%  &    8.01 &   1.6986 &  67.9710 &100\%  &   15.85 &   7.1758   &    67.9710\\ \hline
      80& 77.9747 &100\%  &    8.00 &   3.2806 &  77.9747 &100\%  &   14.80 &  20.4195   &   77.9747\\ \hline
      90& 87.9775  &100\% &    8.01 &   6.1458 &  87.9775  &100\% &   15.09 &  45.8924   &    87.9775    \\ \hline
      100&97.9798  &100\% &    8.00 &  13.7736 &  97.9798  &100\% &   15.42 &  89.7867   &    97.9798  \\ \hline
\end{tabular}
 \label{HGanacontablecomplete}
\end{table}


From Table \ref{HGanacontablecomplete}, we can see that {\tt FTR} takes less iterations and hence less CPU time than that of {\tt SNOPT}, and the upper bound given in \reff{upperbound} is quite tight as it is pretty close to the value from computation. In addition, as the hypergraph $K_n^-$ is well connected by definition, the analytic connectivity is relatively high comparing to all the others in this section, which again verify that the analytic connectivity is a good choice to measure the connectivity of hypergraphs. However, as one can see from Tables \ref{HGanacontable10} and \ref{HGanacontablecomplete}, big analytic connectivities of hypergraphs result in more CPU time for the corresponding hypergraphs with the same $n$.

%

%
%
%
%
%
%

\section{Conclusions}
\label{section:conclusions}

In this paper, we have exploited properties on the analytic connectivity and have shown that several structured uniform hypergraphs attain their analytic connectivities at vertices of the minimum degrees. To efficiently compute the analytic connectivity of any general uniform hypergraph, we have proposed a feasible trust region algorithm with global convergence, and have conducted numerical experiments to shown  the advantages of our algorithm in comparison of other existing ones. All the numerical results have verified that the analytic connectivity is a good choice to measure the connectivity of a hypergraph. Moreover, the efficiency of the proposed algorithm makes the extended version of ``Cheeger inequality" in the setting of uniform hypergraphs practically feasible to efficient bound the Cheeger numbers of uniform hypergraphs.


\section*{Acknowledgements}
\label{section:Acknowledgement}

 The first author thank Professor Wei Li and Miss Lizhu Sun for their
 useful discussions on Conjecture 3.1.
 The first, third and forth authors  were supported in part by the Research Grants Council (RGC) of Hong Kong (Project C1007-15G),
 the second  author's work was supported by the National Natural Science Foundation of China (11301022,11431002),
 and the third author's work was supported by Grants No. PolyU 501212, 501913, 15302114 and 15300715.





\end{document}